\documentclass[10pt, reqno]{amsart}
\usepackage{amssymb}
\usepackage{amsmath,amscd}
\usepackage[initials,nobysame,alphabetic]{amsrefs}
\usepackage{amsfonts}
\usepackage{mathrsfs}
\usepackage{mathpazo}
\usepackage[arrow,matrix,curve,cmtip,ps]{xy}

\usepackage{amsthm}

\usepackage{float}
\usepackage{hyperref}
\usepackage{graphics}
\usepackage{graphicx}
\usepackage{verbatim}
\usepackage{multirow}
\usepackage{tikz}
\usepackage{enumerate}
\usepackage{eufrak}
\allowdisplaybreaks

\newtheorem{theorem}{Theorem}[section]

\newtheorem{corollary}[theorem]{Corollary}

\newtheorem*{theorem*}{Theorem}
\theoremstyle{remark}
\newtheorem{remark}[theorem]{Remark}
\newtheorem{definition}[theorem]{Definition}

\numberwithin{equation}{section}

\newcommand{\R}{\mathbb{R}}

\newcommand{\E}{\mathbb{E}}
\newcommand{\Prob}{\mathbb{P}}

\newcommand{\mytilde}{\raise.17ex\hbox{$\scriptstyle\mathtt{\sim}$}}

\begin{document}
\title[The Principal-Agent Problem]{The Principal-Agent Problem With Time Inconsistent Utility Functions}

\author{Boualem Djehiche and Peter Helgesson}

\address{Department of Mathematics \\ KTH Royal Institute of Technology \\ 100 44, Stockholm \\ Sweden}
\email{boualem@kth.se}
\address{Department of Mathematics \\ Chalmers University of Technology \\ 412 96, G\"{o}teborg \\ Sweden}
\email{helgessp@chalmers.se}


\date{\today}



\begin{abstract}

    In this paper we study a generalization of the continuous time Principal-Agent problem allowing for time inconsistent utility functions, for instance of mean-variance type. Using recent results on the Pontryagin maximum principle for FBSDEs we suggest a method of characterizing optimal contracts for such models. To illustrate this we consider a fully solved explicit example in the linear quadratic setting.

\end{abstract}

\maketitle

\tableofcontents

\section{Introduction}
Risk management or the problem of finding an optimal balance between expected returns and risk taking is a central topic of research within finance and economics. Applications such as portfolio optimization, optimal stopping and liquidation problems have been of particular interest in the literature. In such applications it is common to consider utility functions of mean-variance type. Mean-variance utility functions constitute an important subclass of the so called time inconsistent utility functions for which the Bellman principle of dynamic programming does not hold. Problems involving such utilities can therefore not be approached by the classical Hamilton-Jacobi-Bellman equation. This question has been addressed in Andersson and Djehiche \cite{MR2784835}, Bj\"{o}rk and Murgoci \cite{timeincon}, Bj\"{o}rk, Murgoci and Zhou \cite{MR3157686}, Djehiche and Huang \cite{subgamperfect}, Ekeland and Lazrak \cite{Ekeland} and Ekeland and Pirvu \cite{MR2461340}. In this paper we develop a method of studying a mean-variance setting of the celebrated Principal-Agent problem by means of the stochastic generalization of Pontryagin's maximum principle.\\
\indent The precise structure of the Principal-Agent problem goes as follows: The principal employs an agent to manage a certain well-defined noisy asset over a fixed period of time. In return for his/her effort the agent receives a compensation according to some agreement, set before the period starts. It could for instance involve a lump-sum payment at the end of the period, a continuously paying cash-flow during the period, or both. Depending on what information the principle has at hand to form an agreement, one distinguishes between two cases; the \textit{Full Information}- and the \textit{Hidden Action-problem}. The full information case differs from the hidden action case in that the principal can observe the actions of the agent in addition to the evolution of the asset. Therefore, under full information the principal is allowed to tailor a contract based on both outcome and effort, not only outcome as for hidden actions. In both cases the contract is constrained by the agent via a so called \textit{participation constraint}, clarifying the minimum requirements of the agent to engage in the project. Under hidden action the contract is further constrained by the \textit{incentive compatibility} condition, meaning that as soon as a contract is assigned the agent will act as to maximize his/her own utility and not necessarily that of the principal.\\
\indent The pioneering paper in which the Principal-Agent problem first appears is Holmstr\"{o}m and Milgrom \cite{MR882097}. They study a time continuous model over a finite period in which the principle and the agent both optimize exponential utility functions. The principal rewards the agent at the end of the period by a lump-sum payment. As a result they find that the optimal contract is linear with respect to output. The paper \cite{MR882097} is generalized in \cite{MR1252335} by Sh\"{a}ttler and Sung to a mathematical framework that uses methods from dynamic programming and martingale theory to characterize contract optimality.\\
\indent The interest in continuous time models of the Principal-Agent problem has grown substantially since the first studies appeared. In \cite{MR2465709}, \cite{MR2433118}, \cite{Westerfield}, \cite{Williams} (only to mention a few) the authors analyze continuous time models in a classical setting, i.e. having one principal and one agent. Such models are also covered in the recent book \cite{MR2963805} by Cvitani{\'c} and Zhang. Other models such as various multiplayer versions have been studied for instance in \cite{Kang} and \cite{MR2454676}.\\
\indent Our goal is to characterize optimal contracts in the classical setting of Principal-Agent problem under hidden action for time inconsistent utility functions. We consider two different modeling possibilities; \textit{Hidden Action in the weak formulation} and \textit{Hidden Contract in the strong formulation}. In the first model the agent has full information of the mechanisms behind the cash-flow and the principal wishes to minimize his/her mean-variance utility. In the latter model the agent does not know the structure of the cash-flow and has to protect him/her-self from high levels of risk by an additional participation constraint of variance type. To the best of our knowledge this has not previously been addressed in the literature. In order to carry the program through we use recent generalizations of Pontryagin's stochastic maximum principle. The idea is to consider the Principal-Agent problem as a sequential optimization problem. We first consider the Agent's problem of characterizing optimal choice of effort. Then we proceed to the Principal's problem which, by incentive compatibility, becomes a constrained optimal control problem of a forward-backward stochastic differential equation (from now on FBSDE). A similar scheme was considered in \cite{PA} but without the non-standard mean-variance consideration. Optimal control with respect to mean-variance utility functions have previously been studied in for instance Zhou and Li \cite{MR1751306} and Andersson and Djehiche \cite{MR2784835}.\\
\indent 
\indent In the present literature of the Principal-Agent problem the paper closest to ours is \cite{Williams}, in which a similar maximum principle approach is used. The setting is classical (without time inconsistent utility functions) and the author finds a characterization for optimal choice of effort in the Agent's problem. The full model involving the constrained Principal's problem, however, is not considered. Our contribution to the existing literature should be regarded as mathematical rather than economical. We present a general framework for solving a class of Principal-Agent problems, without claiming or investigating possible economical consequences. The main results of our study are presented in Theorem \ref{HiddenActionMainThm} and Theorem \ref{HiddenContractMainThm} in which a full characterization of optimal contracts is stated for two different models.\\
\indent The paper is organized as follows: In Section 2 we introduce the mathematical machinery from stochastic optimal control theory that is necessary for our purposes. Mean-variance maximum principles are then derived in Section 3 by results from Section 2 in two different but related cases. Section 4 is devoted to fit the methods from previous sections into a Principle-Agent framework. We consider two different models under Hidden Action and find necessary conditions for optimality. Finally in Section 5 we make the general scheme of Section 4 concrete by a simple and fully solved example in the linear-quadratic (LQ)-setting.

\section{Preliminaries}
\label{Preliminaries}

Let $T > 0$ be a fixed time horizon and $(\Omega, \mathcal{F}, \mathbb{F}, \mathbb{P})$ be a filtered probability space satisfying the usual conditions on which a $1$-dimensional Brownian motion $W = \{ W_t \}_{t \geq 0}$ is defined. We let $\mathbb{F}$ be the natural filtration generated by $W$ augmented by all $\mathbb{P}$-null sets $\mathcal{N}_{\mathbb{P}}$, i.e.  $\mathbb{F} = \mathcal{F}_t \bigvee \mathcal{N}_{\mathbb{P}}$ where $\mathcal{F}_t := \sigma(\{ W_s \}: 0 \leq s \leq t)$.\\
\indent Consider the following control system of mean-field type:
\begin{equation*}
\left\{
    \begin{array}{l}
    dx(t) = b(t,x(t),\mathbb{E}[x(t)],s(t)) dt + \sigma (t,x(t),\mathbb{E}[x(t)]) dW_t, \hspace{0.5cm} t \in (0,T]\\
    x(0) = x_0
    \end{array}
\right.
\end{equation*}
with a cost functional of the form
\begin{equation}
    \mathcal{J}(s(\cdot)) := \mathbb{E} \left[ \int_0^T f(t,x(t),\mathbb{E}[x(t)],s(t)) dt + h(x(T),\mathbb{E}[x(T)]) \right],
    \label{IntroCostfunc1}
\end{equation}
where $b:[0,T] \times \R \times \R \times S \rightarrow \R$, $\sigma : [0,T] \times \R\times \R \rightarrow \R$, $f: [0,T] \times \R \times \R \times S \rightarrow \R$ and $h: \R \times \R \rightarrow \R $ and $S \subset \R$ is a non-empty subset. The control $s(\cdot)$ is \textit{admissible} if it is an $\mathbb{F}$-adapted and square-integrable process taking values in $S$. We denote the set of all such admissible controls by $\mathcal{S}[0,T]$. In order to avoid technicalities in regularity that are irrelevant for our purposes we state the following assupmtion.\\
\newline
\textbf{Assumption 1.} The functions $b$, $\sigma$, $f$ and $h$ are  $C^1$ with respect to $x$ and $\tilde{x}$, where $\tilde{x}$ denotes the explicit dependence of $\E[x(\cdot)]$. Moreover, $b$, $\sigma$, $f$ and $h$ and their first order derivatives with respect to $x$ and $\tilde{x}$ are bounded and continuous in $x$, $\tilde{x}$ and $s$.\vspace{0.1cm}\\
\noindent We are interested in the following optimal control problem:\\
\newline
\textbf{Problem (S).} Minimize (\ref{IntroCostfunc1}) over $\mathcal{S}[0,T]$.\\
\newline
Any $\bar{s}(\cdot) \in \mathcal{S}[0,T]$ satisfying
\begin{equation*}
    \mathcal{J}(\bar{s}(\cdot)) = \smash{\displaystyle\inf_{s(\cdot) \in \mathcal{S}[0,T]}} \mathcal{J}(s(\cdot)) \vspace{0.2cm}
\end{equation*}
is called an \textit{optimal control} and the corresponding $\bar{x}(\cdot)$ is called the \textit{optimal state process}. We will refer to $(\bar{x}(\cdot), \bar{s}(\cdot))$ as an \textit{optimal pair}.\vspace{0.1cm}\\
\indent The following stochastic maximum principle for characterizing optimal pairs in problem (S) was found by Buckdahn, Djehiche and Li in \cite{MR2822408}.
\begin{theorem}[The Stochastic Maximum Principle]
    Let the conditions in Assumption 1 hold and consider an optimal pair $(\bar{x}(\cdot), \bar{s}(\cdot))$ of problem (S). Then there exists a pair of processes $(p(\cdot), q(\cdot)) \in L_{\mathcal{F}}^2(0,T;\R) \times (L_{\mathcal{F}}^2(0,T;\R))$ satisfying the adjoint equation
    \begin{equation}
    \left\{
    \begin{array}{lll}
        dp(t) = \vspace{0.1cm}-\left\{ b_x(t, \bar{x}(t), \E[\bar{x}(t)], \bar{s}(t)) p(t) + \E[b_{\tilde{x}}(t, \bar{x}(t), \E[\bar{x}(t)], \bar{s}(t)) p(t)]\right. \\
	\vspace{0.1cm}\hspace{1.2cm}\left. + \sigma_x (t,x(t),\mathbb{E}[\bar{x}(t)]) q(t) + \E[\sigma_{\tilde{x}} (t,x(t),\mathbb{E}[\bar{x}(t)]) q(t)]\right.\\
	\vspace{0.1cm}\hspace{1.2cm}\left. - f_x(t, \bar{x}(t), \mathbb{E}[\bar{x}(t)]), \bar{s}(t)) - \E[f_{\tilde{x}}(t, \bar{x}(t), \mathbb{E}[\bar{x}(t)]), \bar{s}(t))] \right\} dt\\
	\vspace{0.1cm}\hspace{1.2cm}+ q(t) dW_t, \\
        p(T) = -h_x(\bar{x}(T), \E[\bar{x}(T)]) - \E[h_{\tilde{x}}(\bar{x}(T), \E[\bar{x}(T)])],
    \end{array}
    \right.
    \end{equation}
    such that
    \begin{equation}
        \bar{s}(t)= \smash{\displaystyle\arg\max_{s \in S}}\hspace{0.1cm}\mathcal{H}(t, \bar{x}(t), s, p(t), q(t)), \hspace{0.3cm} \text{a.e. } t \in [0,T], \hspace{0.3cm} \Prob\text{-a.s.}
    \end{equation}
    where the Hamiltonian function $\mathcal{H}$ is given by
    \begin{equation}
        \mathcal{H} (t, x, s, p, q) := b(t,x,\E[x],s) \cdot p  + \sigma(t,x,\E[x]) \cdot q - f(t,x,\E[x],s)
    \end{equation}
    for $(t,x,s,p,q) \in [0,T] \times \R \times S \times \R \times \R$.
    \label{PrelThmSMP1}
\end{theorem}
\begin{remark}
It is important to remember that Theorem \ref{PrelThmSMP1} merely states a set of necessary conditions for optimality in (S). It does not claim the existence of an optimal control. Existence theory of stochastic optimal controls (both in the strong and the weak sense) has been a subject of study since the sixties (see e.g. \cite{MR0192946}) and, at least in the case of strong solutions, the results seem to depend a lot upon the statement of the problem. In the weak sense an account of existence results is  to be found in \cite{MR1696772} (Theorem 5.3, p. 71).
\label{PrelRmkExistence}
\end{remark}
\begin{remark}
    Restricting the space $U$ to be convex allows for a diffusion coefficient of the form $\sigma(t,x, \E[x],s)$, without changing the conclusion of Theorem \ref{PrelThmSMP1}. In the case of a non-convex control space the stochastic maximum principle with controlled diffusion was proven in \cite{MR1051633} and requires the solution of an additional adjoint BSDE. We choose to leave this most general maximum principle as reference in order to keep the presentation clear.
    \label{PrelRemarkThmSMP1}
\end{remark}
As pointed out in Remark \ref{PrelRmkExistence} it is a non-trivial task to prove the existence of an optimal pair $(\bar{x}(\cdot), \bar{s}(\cdot))$ in a general stochastic control model. Under the additional assumptions;\\
\newline
\textbf{Assumption 2.} The control domain $S$ is a convex body in $\R$. The maps $b$, $\sigma$ and $f$ are locally Lipschitz in $u$ and their derivatives in $x$ and $\tilde{x}$ are continuous in $x$, $\tilde{x}$ and $s$,\\
\newline
the following theorem provides sufficient conditions for optimality in (S).
\begin{theorem}[Sufficient Conditions for Optimality]
Under Assupmptions 1 and 2 let $(\bar{x}(\cdot), \bar{s}(\cdot), p(\cdot), q(\cdot))$ be an admissible 4-tuple. Suppose that $h$ is convex and further that $\mathcal{H} (t, \cdot, \cdot, \cdot, p(t), q(t))$ is concave for all $t \in [0,T]$ $\mathbb{P}$-a.s. and
\begin{equation*}
    \bar{s}(t)= \smash{\displaystyle\arg\max_{s \in S}}\hspace{0.1cm}\mathcal{H}(t, \bar{x}(t), \E[\bar{x}(t)], s, p(t), q(t)), \hspace{0.3cm} \text{a.e. } t \in [0,T], \hspace{0.3cm} \Prob\text{-a.s.}
\end{equation*}
Then $(\bar{x}(\cdot), \bar{s}(\cdot))$ is an optimal pair for problem (S).
\end{theorem}
\indent The stochastic maximum principle has since the early days of the subject (in pioneering papers by e.g. Bismut \cite{MR0469466} and Bensoussan \cite{MR705931}) developed a lot and does by now apply to a wide range of problems more general than (S) (see for instance \cite{MR1051633}, \cite{MR2784835} \cite{MR2822408}, \cite{risksensitive}). For our purposes we need a refined version of Theorem \ref{PrelThmSMP1}, characterizing optimal controls in a FBSDE-dynamical setting under state constraints. More precisely we wish to consider a stochastic control system of the form
\begin{equation}
    \left\{
    \begin{array}{l}
         \vspace{0.1cm}dx(t) = b(t, \Theta(t), s(t))dt + \sigma(t, \Theta(t))dW_t\\
         \vspace{0.1cm}dy(t) = -c(t, \Theta(t), s(t))dt + z(t) dW_t\\
         x(0) = x_0, \hspace{0.2cm} y(T) = \varphi(x(T)),
    \end{array}
    \right.
    \label{PrelFBSDEsyst}
\end{equation}
where $b, \sigma , c: [0,T] \times \R^6 \times S \rightarrow \R$ and $\varphi : \R \rightarrow \R$, with respect to a cost-functional of the form
\begin{equation}
	\begin{array}{l}
		\displaystyle\vspace{0.1cm}\mathcal{J}(s(\cdot)) := \mathbb{E} \left[ \int_0^T f(t,\Theta(t),s(t)) dt + h(x(T), \E[x(T)]) + g(y(0)) \right],
	\end{array}
    \label{PrelCostfunc2}
\end{equation}
and a set of state constraints
\begin{equation}
    \begin{array}{l}
    \displaystyle\E \left[ \int_0^T \mathbf{F}(t, \Theta(t), s(t)) dt + \mathbf{H}(x(T), \E[x(T)]) + \mathbf{G}(y(0)) \right] :=\vspace{0.3cm}\\
    \left(
    \begin{array}{l}
        \displaystyle\E \left[ \int_0^T f^1(t, \Theta(t), s(t)) dt + h^1(x(T), \E[x(T)]) + g^1(y(0)) \right]\\
        \hspace{4cm}\vdots\\
        \displaystyle\E \left[ \int_0^T f^l(t, \Theta(t), s(t)) dt + h^l(x(T), \E[x(T)]) + g^l(y(0)) \right]
    \end{array}
    \right) \in \Lambda,
    \end{array}
    \label{PrelStateConstraints}
\end{equation}
for some closed and convex set $\Lambda \subseteq \R^l$. In the above expressions we have introduced
\begin{equation*}
	\Theta(t) := (x(t), y(t), z(t), \E[x(t)], \E[y(t)], \E[z(t)]),
\end{equation*}
in order to avoid unnecessarily heavy notation. The optimal control problem is:\\
\newline
\textbf{Problem (SC).} Minimize (\ref{PrelCostfunc2}) subject to the state constraints (\ref{PrelStateConstraints}) over the set $\mathcal{S}[0,T]$.\vspace{0.1cm}\\
\newline
To get a good maximum principle for (SC) we require some further regularity conditions ensuring solvability of (\ref{PrelFBSDEsyst}). These conditions are listed in the following assumptions and can be found in \cite{MR3178298}.\\
\newline
\textbf{Assumption 3.} The functions $b, \sigma, c$ are continuously differentiable and Lipschitz continuous in $\Theta$, the functions $h$, $g$, $h^i$, $g^i$ are continuously differentiable in $x$ and $y$ respectively, and they are bounded by $C(1 + |x| + |y| + |z| + |\tilde{x}| + |\tilde{y}| + |\tilde{z}| + |s|)$, $C(1 + |x|)$ and $C(1 + |y|)$ respectively.\\
\newline
\textbf{Assumption 4.} All derivatives in Assumption 4 are Lipschitz continuous and bounded.\\
\newline
\textbf{Assumption 5.} For all $\Theta \in \R^6$, $s \in \mathcal{S}$, $A(\cdot , \Theta, s) \in L_{\mathcal{F}}^2(0,T;\R^3)$, where we have $A(t,\Theta, s) :=(c(t,\Theta, s), b(t,\Theta, s), \sigma(t,\Theta))$ and
\begin{equation*}
	L_{\mathcal{F}}^2(0,T;\R^k) := \left\{ \psi : [0,T] \times \Omega \rightarrow \R^k \Big| \hspace{0.1cm} \psi \hspace{0.1cm} \text{ is } \mathbb{F}\text{-adapted and } \E \left[ \int_0^T |\psi|^2 dt \right] < \infty \right\},
\end{equation*}
and for each $x \in \R$, $\varphi (x) \in L^2_{\mathcal{F}} (\Omega ; \R)$. Furthermore, there exists a constant $C > 0$ such that
\begin{equation*}
	\left\{
	\begin{array}{l}
		\vspace{0.1cm}|A(t,\Theta_1, s) - A(t,\Theta_2, s)| \leq C|\Theta_1 - \Theta_2|, \hspace{0.2cm} \Prob\text{-a.s. and for a.e. } t \in [0,T],\\
		\vspace{0.1cm}|\varphi(x_1) - \varphi(x_2)| \leq C|x_1 - x_2|, \hspace{0.2cm} \Prob\text{-a.s},\\
		\text{for all } \Theta_1, \Theta_2 \in \R^6.
	\end{array}
	\right.
\end{equation*}
\newline
\textbf{Assumption 6.} The functions $A$ and $\varphi$ satisfy the following monotonicity conditions:
\begin{equation*}
\left\{
	\begin{array}{l}
		\E \langle A(t,\Theta_1, s) - A(t,\Theta_2, s), \Theta_1 - \Theta_2 \rangle \leq \beta \E |\Theta_1 - \Theta_2|^2, \hspace{0.2cm} \Prob \text{-a.s}\vspace{0.1cm}\\
		\langle \varphi (x_1) - \varphi(x_2), x_1 - x_2 \rangle \geq \mu |x_1 - x_2|^2
	\end{array}
\right.
\end{equation*}
for all $\Theta_1, \Theta_2 \in \R^6$, $x_1, x_2 \in \R$\\
\newline
In the spirit of \cite{MR3178298} we are now ready to formulate the state constrained stochastic maximum principle for fully coupled FBSDEs of mean-field type.
\begin{theorem}[The State Constrained Maximum Principle]
    Let Assumptions 3-6 hold and assume $\Lambda \subseteq \R^l$ to be a closed and convex set. If $(\bar{x}(\cdot), \bar{y}(\cdot), \bar{z}(\cdot), \bar{s}(\cdot))$ is an optimal 4-tuple of problem (SC), then there exists a vector $(\lambda_0, \lambda) \in \R^{1+l}$ such that
    \begin{equation}
        \lambda_0 \geq 0, \hspace{0.5cm} |\lambda_0|^2 + |\lambda|^2 = 1,
        \label{PrelLagMult}
    \end{equation}
    satisfying the transversality condition
    \begin{equation}
        \langle \lambda , v - \E \left[ \int_0^T \mathbf{F}(t, {\bar x}(t), {\bar y}(t), {\bar z}(t), {\bar s}(t)) dt + \mathbf{H}({\bar x}(T)) + \mathbf{G}({\bar y}(0)) \right] \rangle \geq 0, \hspace{0.3cm} \forall v \in \Lambda
        \label{PrelTransversalityCond}
    \end{equation}
    and a 3-tuple $(r(\cdot), p(\cdot), q(\cdot)) \in L_{\mathcal{F}}^2(\Omega; C([0,T]; \mathbb{R})) \times L_{\mathcal{F}}^2(\Omega; C([0,T]; \mathbb{R})) \times L_{\mathcal{F}}^2(0,T;\mathbb{R})$ of solutions to the adjoint FBSDE
    \begin{equation*}
        \begin{array}{l}
        dr(t)= \left\{ c_y(t) r(t) - b_y(t) p(t) - \sigma_y(t) q(t) + \sum_{i=0}^l \lambda_i f_y^i(t)\right.\\
	\vspace{0.2cm}\hspace{2cm}\left. + \E [c_{\tilde{y}}(t) r(t) - b_{\tilde{y}}(t) p(t) - \sigma_{\tilde{y}}(t) q(t) + \sum_{i=0}^l \lambda_i f_{\tilde{y}}^i(t)] \right\} dt\\
        \vspace{0.2cm}\hspace{0.9cm}+\left\{ c_z(t) r(t) - b_z (t)p(t) - \sigma_z(t) q(t) + \sum_{i=0}^l \lambda_i f_z^i(t)\right.\\
	\vspace{0.2cm}\hspace{2cm}\left. + \E [c_{\tilde{z}}(t) r(t) - b_{\tilde{z}}(t) p(t) - \sigma_{\tilde{z}}(t) q(t) + \sum_{i=0}^l \lambda_i f_{\tilde{z}}^i(t)] \right\} dW_t,\\
    \end{array}
    \end{equation*}
    \begin{equation}
    \begin{array}{l}
        \vspace{0.2cm}dp(t) = -\left\{ -c_x(t) r(t) + b_x(t) p(t) + \sigma_x(t) q(t) - \sum_{i=0}^l \lambda_i f_x^i(t)\right.\\
	\vspace{0.2cm}\hspace{1.7cm}\left. + \E [-c_{\tilde{x}}(t) r(t) + b_{\tilde{x}}(t) p(t) + \sigma_{\tilde{x}}(t) q(t) - \sum_{i=0}^l \lambda_i f_{\tilde{x}}^i(t)] \right\} dt\\
	\vspace{0.2cm}\hspace{9.5cm}+ q(t) dW_t,\\
        \vspace{0.2cm}r(0) = \sum_{i=0}^l \lambda_i \E [ g^i({\bar y}(0)) ],\\
	p(T) = -\varphi_x({\bar x}(T))r(T) - \sum_{i=0}^l \lambda_i ( h_x^i({\bar x}(T), \E [{\bar x}(T)]) + \E [h_{\tilde{x}}^i({\bar x}(T), \E [{\bar x}(T)])]),
        \end{array}
        \label{PrelFBSDEConstSMP}
    \end{equation}
    such that
    \begin{equation*}
        \bar{s}(t)= \smash{\displaystyle\arg\max_{s \in S}}\hspace{0.1cm} \mathcal{H}(t, \bar{\Theta}(t),s,r(t),p(t),q(t),\lambda_0,\lambda) \hspace{0.3cm} \text{a.e. } t \in [0,T], \hspace{0.3cm} \Prob\text{-a.s.}
    \end{equation*}
    where the Hamiltonian function $\mathcal{H}$ is given by
    \begin{equation*}
        \begin{array}{l}
        \vspace{0.1cm}\mathcal{H}(t,\Theta,s,r,p,q,\lambda_0,\lambda) :=\\
         \hspace{0.5cm} -r \cdot c(t,\Theta,s) + p \cdot b(t,\Theta ,s)  + q \cdot \sigma(t,\Theta) - \sum_{i=0}^l \lambda_i f^i(t,\Theta,s).
        \end{array}
    \end{equation*}
    \label{PrelConstrainedSMPthm}
\end{theorem}
\begin{remark}
As in Remark \ref{PrelRemarkThmSMP1}, analogue principles also hold in Theorem \ref{PrelConstrainedSMPthm}.
\end{remark}
\begin{remark}
The maximum principle in Theorem \ref{PrelConstrainedSMPthm} without state constraints is an easy extension of the same result in \cite{MR3178298} and follows the proof \textit{mutatis mutandis}. Extending the result to allow for state constraints is a standard procedure and can be found for instance in \cite{PA}.
\end{remark}

\section{Utilities of mean-variance type}
\label{MeanVarianceSMP}
We are now going to fit the methods presented in Section 2 to a mean-variance framework, i.e. we want to control a FBSDE of mean-field type (\ref{PrelFBSDEsyst}) with respect to either of the following two cases:\\
\newline
\textbf{(i).} Minimize
\begin{equation}
    \begin{array}{l}
        \vspace{0.2cm}\displaystyle \mathcal{I}(u) := -\E \left[ \int_0^T U (t, \Theta(t), s(t)) dt +V(x(T)) \right]\\
        \hspace{0.5cm}\displaystyle+ \frac{r}{2} \text{Var} \left( \int_0^T \Phi (t, \Theta(t), s(t)) dt + \Psi(x(T)) \right),
    \end{array}
    \label{MVCostFunc1}
\end{equation}
\noindent over $\mathcal{S}[0,T]$ for some risk aversion $r > 0$.\\
\newline
\textbf{(ii).} Minimize
\begin{equation}
    \mathcal{J}(u) := \mathbb{E} \left[ \int_0^T U(t,\Theta(t),s(t)) dt + V(x(T)) \right]
    \label{MVCostFunc3}
\end{equation}
over $\mathcal{S}[0,T]$ subject to a set of state constraints (compare (\ref{PrelStateConstraints})), including statements of the form
\begin{equation}
    \text{Var} \left( \int_0^T \Phi(t,\Theta(t),s(t)) dt + \Psi(x(T)) \right) \leq R_0.
    \label{MVConstraintVar}
\end{equation}
In order to carry this through we introduce the auxiliary process
\begin{equation*}
    \eta (t) := \int_0^t \Phi (\tau, \Theta(\tau), s(\tau)) d\tau + \Psi(x(t)),
\end{equation*}
which by It\^{o}'s Lemma solves the SDE
\begin{equation}
\left\{
    \begin{array}{l}
        d \eta (t) = \left\{ \Phi(t) + b(t) \cdot \Psi'(x(t)) + \frac{\sigma^2(t)}{2} \cdot \Psi''(x(t)) \right\} dt + \sigma(t) \cdot \Psi'(x(t)) dW_t,\\
        \eta(0) = 0.
    \end{array}
    \label{MVAux}
\right.
\end{equation}
Here we adopt the simpler notational convention $\Phi(t) := \Phi (t, \Theta(t), s(t))$. By considering (\ref{PrelFBSDEsyst}) with (\ref{MVAux}) as an augmented dynamics we may rewrite (\ref{MVCostFunc1}) (or analogously for the state constraint (\ref{MVConstraintVar}) )as
\begin{equation}
    \mathcal{I}(s) := \E \left[ -\int_0^T U (t, \Theta(t), s(t)) dt + V(x(T)) + \frac{r}{2}(\eta (T) - \E[\eta (T)])^2 \right].
    \label{MVCostFunc2}
\end{equation}
An optimal control problem involving the cost functional (\ref{MVCostFunc2}) is within the framework of Section 2, in particular Theorem \ref{PrelFBSDEConstSMP}.\\
\indent For the Principal-Agent problem we are interested in the following:\\
\newline
\noindent\textbf{(MV1).} Minimize (\ref{MVCostFunc1}) subject to the state constraint
\begin{equation*}
    \E \left[ \int_0^T u (t, \Theta(t), s(t)) dt + v(x(T)) \right] \leq W_0,
\end{equation*}
for some finite $W_0 \in \R$ over $\mathcal{S}[0,T]$.\\
\newline
\noindent\textbf{(MV2).} Minimize (\ref{MVCostFunc3}) subject to the state constraints
\begin{equation*}
    \left\{
    \begin{array}{l}
        \displaystyle\vspace{0.1cm}\mathcal{J}_{E}(s) := \E \left[ \int_0^T u (t, \Theta(t), s(t)) dt + v(x(T)) \right] \leq W_0,\\
        \displaystyle\mathcal{J}_{V}(s) := \text{Var}\left( \int_0^T \phi (t, \Theta(t), s(t)) dt + \psi(x(T)) \right) \leq R_0,
    \end{array}
    \right.
\end{equation*}
for some finite $W_0 < 0$ and $R_0 > 0$ over $\mathcal{S}[0,T]$.\\
\newline
It is now an easy task to formulate the stochastic maximum principles that characterize optimality in (MV1) and (MV2) respectively. In the two Corollaries that follow we adopt the vector notation:
\begin{equation*}
    \textbf{B}(t) := \begin{pmatrix} b(t) \\ \Phi(t) + b(t) \cdot \Psi'(x(t)) + \frac{\sigma^2(t)}{2} \cdot \Psi''(x(t)) \end{pmatrix}, \hspace{0.2cm} \boldsymbol\Sigma(t) := \sigma(t) \begin{pmatrix} 1 \\ \Psi'(t) \end{pmatrix},
\end{equation*}
\begin{equation*}
    \textbf{b}(t) :=  \begin{pmatrix} b(t) \\ \phi(t) + b(t) \cdot \psi'(x(t)) + \frac{\sigma^2(t)}{2} \cdot \psi''(x(t)) \end{pmatrix}, \hspace{0.2cm} \boldsymbol\sigma(t) := \sigma(t) \begin{pmatrix} 1 \\ \psi'(t) \end{pmatrix},
\end{equation*}
and
\begin{equation*}
         \textbf{x}(t) := \begin{pmatrix} x(t) \\ \eta(t) \end{pmatrix}.
\end{equation*}
\begin{corollary}[The Stochastic Maximum Principle for MV1]
Let Assumptions 3-6 hold and let $\Psi(\cdot)$ be three times differentiable. If $(\bar{x}(\cdot), \bar{y}(\cdot), \bar{z}(\cdot), \bar{s}(\cdot))$ is an optimal 4-tuple of (MV1), then there exists a vector $(\lambda_A, \lambda_P) \in \R^{2}$ such that
    \begin{equation}
        \lambda_P \geq 0, \hspace{0.5cm} \lambda_P^2 + \lambda_A^2 = 1,
        \label{PrelLagMult}
    \end{equation}
and a 3-tuple $(R(\cdot), P(\cdot), Q(\cdot)) \in L_{\mathcal{F}}^2(\Omega; C([0,T]; \mathbb{R})) \times (L_{\mathcal{F}}^2(\Omega; C([0,T]; \mathbb{R})))^2 \times (L_{\mathcal{F}}^2(0,T;\mathbb{R}))^2$ of solutions to the adjoint FBSDE
\begin{equation*}
    \begin{array}{l}
        \vspace{0.1cm}dR(t) = \left\{ c_y(t) R(t) - \textbf{B}^T_y(t) \cdot P(t) - \boldsymbol\Sigma^T_y(t) \cdot Q(t) + \lambda_A u_y(t) - \lambda_P U_y(t) \right.\\
        \vspace{0.1cm}\hspace{2cm}+\E[ c_{\tilde{y}}(t) R(t) - \textbf{B}^T_{\tilde{y}}(t) \cdot P(t) - \boldsymbol\Sigma^T_{\tilde{y}}(t) \cdot Q(t) + \lambda_A u_{\tilde{y}}(t) - \lambda_P U_{\tilde{y}}(t) ] \Big\} dt\\
        \vspace{0.1cm}+\left\{ c_z(t) R(t) - \textbf{B}^T_z(t) \cdot P(t) - \boldsymbol\Sigma^T_z(t) \cdot Q(t) + \lambda_A u_z(t) - \lambda_P U_z(t) \right.\\
        \vspace{0.1cm}\hspace{2cm}+\E[ c_{\tilde{z}}(t) R(t) - \textbf{B}^T_{\tilde{z}}(t) \cdot P(t) - \boldsymbol\Sigma^T_{\tilde{z}}(t) \cdot Q(t) + \lambda_A u_{\tilde{z}}(t) - \lambda_P U_{\tilde{z}}(t) ] \Big\} dW_t,\\
    \end{array}
\end{equation*}
\begin{equation}
    \begin{array}{l}
        \vspace{0.1cm}dP(t) = \left\{ c_{\textbf{x}}(t) R(t) - \textbf{B}^T_{\textbf{x}}(t) \cdot P(t) - \boldsymbol\Sigma^T_{\textbf{x}}(t) \cdot Q(t) + \lambda_A u_{\textbf{x}}(t) - \lambda_P U_{\textbf{x}}(t) \right.\\
        \vspace{0.1cm}\hspace{0.1cm}+\E[ c_{\tilde{{\textbf{x}}}}(t) R(t) - \textbf{B}^T_{\tilde{{\textbf{x}}}}(t) \cdot P(t) - \boldsymbol\Sigma^T_{\tilde{{\textbf{x}}}}(t) \cdot Q(t) + \lambda_A u_{\tilde{{\textbf{x}}}}(t) -\lambda_P U_{\tilde{{\textbf{x}}}}(t) ] \Big\} dt + q(t)dW_t,
    \end{array}
    \label{MVCor1FBSDE}
\end{equation}
where
\begin{equation*}
    R(0) = 0, \hspace{0.2cm} P(T) = - \begin{pmatrix} \varphi_x(x(T)) \\ 0 \end{pmatrix} \cdot R(T) - \begin{pmatrix} \lambda_A v'(x(T)) - \lambda_P V'(x(T)) \\ r(\eta(T) - \E[\eta(T)]) \end{pmatrix}
\end{equation*}
such that
\begin{equation*}
\vspace{0.2cm}\bar{s}(t)= \smash{\displaystyle\arg\max_{s \in S}}\hspace{0.1cm} \mathcal{H}(t, \bar{\Theta}(t),s,R(t),P(t),Q(t),\lambda_A,\lambda_P) \hspace{0.3cm} \text{a.e. } t \in [0,T], \hspace{0.3cm} \Prob\text{-a.s.}
\end{equation*}
where the Hamiltonian function $\mathcal{H}$ is given by
\begin{equation*}
    \begin{array}{l}
        \vspace{0.1cm}\mathcal{H}(t,\Theta,s,R,P,Q,\lambda_A,\lambda_P) :=\\
         \hspace{0.2cm} -c(t,\Theta,s) \cdot R + \textbf{B}^T(t,\Theta,s) \cdot P  + \boldsymbol\Sigma^T(t,\Theta,s) \cdot Q - \lambda_A u(t,\Theta,s) + \lambda_P U(t,\Theta,s).
    \end{array}
\end{equation*}
\label{MVSMPMV1}
\end{corollary}
\begin{corollary}[The Stochastic Maximum Principle for MV2]
Let Assumptions 3-6 hold and let $\psi(\cdot)$ be three times differentiable. If $(\bar{x}(\cdot), \bar{y}(\cdot), \bar{z}(\cdot), \bar{s}(\cdot))$ is an optimal 4-tuple of (MV2), then there exists a vector $(\lambda_P, \lambda_E, \lambda_V) \in \R^{3}$ such that
    \begin{equation}
        \lambda_P \geq 0, \hspace{0.5cm} \lambda_P^2 + \lambda_E^2 + \lambda_V^2 = 1,
        \label{PrelLagMult}
    \end{equation}
satisfying the transversality condition
    \begin{equation}
        \begin{pmatrix} \lambda_E \\ \lambda_V \end{pmatrix}^T \cdot \begin{pmatrix} \vspace{0.1cm} v_1 - \E \left[ \int_0^T u (t, \Theta(t), s(t)) dt + v(x(T)) \right] \\ v_2 - \text{Var}\left( \int_0^T \phi (t, \Theta(t), s(t)) dt + \psi(x(T)) \right) \end{pmatrix} \geq 0, \hspace{0.3cm} \forall \hspace{0.1cm} v_1 \leq W_0, v_2 \leq R_0
        \label{PrelTransversalityCond}
    \end{equation}
and a 3-tuple $(R(\cdot), P(\cdot), Q(\cdot)) \in L_{\mathcal{F}}^2(\Omega; C([0,T]; \mathbb{R})) \times (L_{\mathcal{F}}^2(\Omega; C([0,T]; \mathbb{R})))^2 \times (L_{\mathcal{F}}^2(0,T;\mathbb{R}))^2$ of solutions to the adjoint FBSDE
\begin{equation*}
    \begin{array}{l}
        \vspace{0.1cm}dR(t) = \left\{ c_y(t) R(t) - \textbf{b}^T_y(t) \cdot P(t) - \boldsymbol\sigma^T_y(t) \cdot Q(t) + \lambda_E u_y(t) + \lambda_P U_y(t) \right.\\
        \vspace{0.1cm}\hspace{2cm}+\E[ c_{\tilde{y}}(t) R(t) - \textbf{b}^T_{\tilde{y}}(t) \cdot P(t) - \boldsymbol\sigma^T_{\tilde{y}}(t) \cdot Q(t) + \lambda_E u_{\tilde{y}}(t) + \lambda_P U_{\tilde{y}}(t) ] \Big\} dt\\
        \vspace{0.1cm}+\left\{ c_z(t) R(t) - \textbf{b}^T_z(t) \cdot P(t) - \boldsymbol\sigma^T_z(t) \cdot Q(t) + \lambda_E u_z(t) + \lambda_P U_z(t) \right.\\
        \vspace{0.1cm}\hspace{2cm}+\E[ c_{\tilde{z}}(t) R(t) - \textbf{b}^T_{\tilde{z}}(t) \cdot P(t) - \boldsymbol\sigma^T_{\tilde{z}}(t) \cdot Q(t) + \lambda_E u_{\tilde{z}}(t) + \lambda_P U_{\tilde{z}}(t) ] \Big\} dW_t,\\
    \end{array}
\end{equation*}
\begin{equation}
    \begin{array}{l}
        \vspace{0.1cm}dP(t) = \left\{ c_{\textbf{x}}(t) R(t) - \textbf{b}^T_{\textbf{x}}(t) \cdot P(t) - \boldsymbol\sigma^T_{\textbf{x}}(t) \cdot Q(t) + \lambda_E u_{\textbf{x}}(t) + \lambda_P U_{\textbf{x}}(t) \right.\\
        \vspace{0.1cm}\hspace{0.1cm}+\E[ c_{\tilde{{\textbf{x}}}}(t) R(t) - \textbf{b}^T_{\tilde{{\textbf{x}}}}(t) \cdot P(t) - \boldsymbol\sigma^T_{\tilde{{\textbf{x}}}}(t) \cdot Q(t) + \lambda_E u_{\tilde{{\textbf{x}}}}(t) + \lambda_P U_{\tilde{{\textbf{x}}}}(t) ] \Big\} dt + Q(t)dW_t,
    \end{array}
    \label{MVCor2FBSDE}
\end{equation}
where
\begin{equation*}
    R(0) = 0, \hspace{0.2cm} P(T) = - \begin{pmatrix} \varphi_x(x(T)) \\ 0 \end{pmatrix} \cdot R(T) - \begin{pmatrix} \lambda_P V'(x(T)) + \lambda_E v'(x(T)) \\ 2\lambda_V (\eta(T) - \E[\eta(T)]) \end{pmatrix}
\end{equation*}
such that
\begin{equation*}
\vspace{0.2cm}\bar{s}(t)= \smash{\displaystyle\arg\max_{s \in S}}\hspace{0.1cm} \mathcal{H}(t, \bar{\Theta}(t),s,R(t),P(t),Q(t),\lambda_E,\lambda_P) \hspace{0.3cm} \text{a.e. } t \in [0,T], \hspace{0.3cm} \Prob\text{-a.s.}
\end{equation*}
where the Hamiltonian function $\mathcal{H}$ is given by
\begin{equation*}
    \begin{array}{l}
        \vspace{0.1cm}\mathcal{H}(t,\Theta,s,R,P,Q,\lambda_A,\lambda_P) :=\\
         \hspace{0.2cm} -c(t,\Theta,s) \cdot R + \textbf{b}^T(t,\Theta,s) \cdot P  + \boldsymbol\sigma^T(t,\Theta,s) \cdot Q - \lambda_E u(t,\Theta,s) - \lambda_P U(t,\Theta,s).
    \end{array}
\end{equation*}
\label{MVSMPMV2}
\end{corollary}
\noindent The transversality condition (\ref{PrelTransversalityCond}) specifies which multipliers $(\lambda_E, \lambda_V)$ satisfying (\ref{PrelLagMult}) (given $\lambda_P \in [0,1]$) that are of interest for characterizing optimality in (MV2). If we let
\begin{equation*}
    \vspace{0.1cm}\Lambda_{MV2} := \{ (x,y) \in \R^2 : x \leq W_0, 0 \leq y \leq R_0 \},
\end{equation*}
it is clear by (\ref{PrelTransversalityCond}) that $(\mathcal{J}_E(\bar{s}), \mathcal{J}_V(\bar{s})) \in \partial \Lambda_{MV2}$. This narrows the set of multipliers to five distinct cases:\\
\newline
\textbf{(i).}\hspace{0.2cm} If $\mathcal{J}_V(\bar{s}) = 0$, then $\lambda_E = 0,\hspace{0.1cm} \lambda_V = \sqrt{1 - \lambda_P^2},\hspace{0.1cm} 0 \leq \lambda_P \leq 1$.\\
\textbf{(ii).}\hspace{0.1cm} If $\mathcal{J}_V(\bar{s}) = R_0$, then $\lambda_E = 0,\hspace{0.1cm} \lambda_V = -\sqrt{1 - \lambda_P^2},\hspace{0.1cm} 0 \leq \lambda_P \leq 1$.\\
\textbf{(iii).} If $\mathcal{J}_E(\bar{s}) = W_0$ and $\mathcal{J}_V(\bar{u}) = 0$, then $\lambda_E = -\sqrt{1 - \lambda_P^2}$cos$\theta$,\\
$\left.\hspace{0.7cm}\lambda_V = -\sqrt{1 - \lambda_P^2}\right.$sin$\theta$ $,\hspace{0.1cm} 0 \leq \lambda_P \leq 1, \theta \in \left[ -\frac{\pi}{2}, 0 \right].$\\
\textbf{(iv).} If $\mathcal{J}_E(\bar{s}) = W_0$ and $\mathcal{J}_V(\bar{u}) = R_0$, then $\lambda_E = -\sqrt{1 - \lambda_P^2}$cos$\theta$,\\
$\left.\hspace{0.7cm}\lambda_V = -\sqrt{1 - \lambda_P^2}\right.$sin$\theta$ $,\hspace{0.1cm} 0 \leq \lambda_P \leq 1, \theta \in \left[0, \frac{\pi}{2} \right].$\\
\textbf{(v).}\hspace{0.1cm} If $\mathcal{J}_E(\bar{s}) = W_0$, then $\lambda_E = -\sqrt{1 - \lambda_P^2},\hspace{0.1cm} \lambda_V = 0 ,\hspace{0.1cm} 0 \leq \lambda_P \leq 1$.\\
\newline
Each of the cases (i)-(v) are illustrated in Figure \ref{MVFig1}.
\begin{center}
\begin{figure}[ht!]
\begin{tikzpicture}

\fill[blue!30!white] (-1,0) rectangle (-6,2.5);

\draw (-1,0) -- (-1,2.5);
\draw (-1,2.5) -- (-6,2.5);

\draw[thick,->] (-6,0) -- (1,0) node[anchor=north west] {$\mathcal{J}_E$};
\draw[thick,->] (0,-0.5) -- (0,3.5) node[anchor=south east] {$\mathcal{J}_V$};

\draw[thick,dashed] (-1,0.7) arc (90:180:0.7cm);
\draw[thick,dashed] (-1.7,2.5) arc (180:270:0.7cm);

\fill[black] (-4,0) circle (1.5pt);
\fill[black] (-1,0) circle (1.5pt);
\fill[black] (-1,1.25) circle (1.5pt);
\fill[black] (-1,2.5) circle (1.5pt);
\fill[black] (-4,2.5) circle (1.5pt);

\node at (-1,-0.3) {$W_0$};
\node at (0.3,2.5) {$R_0$};

\draw[thick] (-0.04,2.5) -- (0.04,2.5);

\draw[thick,->] (-4,0) -- (-4,0.7);
\draw[thick,->] (-1,0) -- (-1.5,0.5);
\draw[thick,->] (-1,1.25) -- (-1.7,1.25);
\draw[thick,->] (-1,2.5) -- (-1.5,2);
\draw[thick,->] (-4,2.5) -- (-4,1.8);

\path (-5,2)  node(a)[]  {}
    (-6.5,3.5) node(b)[]  {$\Lambda_{MV2}$};
 \draw[thick] (a) .. controls +(left:0.5cm) and +(down:1cm) .. (b);

\end{tikzpicture}
\caption{The cases (i)-(v) specifying the multipliers $(\lambda_E, \lambda_V)$ that are illustrated by arrows.}
\label{MVFig1}
\end{figure}
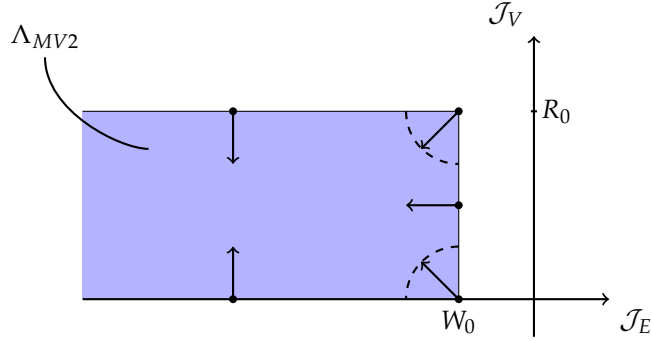
\end{center}

\section{The Principal-Agent Problem}
\label{TheProblem}
We are now ready to state the Principal-Agent problem in the framework of Sections \ref{Preliminaries} and \ref{MeanVarianceSMP} and thereby develop a scheme for characterizing optimality. In the present literature two types of models seem to be the most popular; \textit{The Full Information} case and the \textit{Hidden Action} case. Our treatment will focus on the Hidden Action regime although similar techniques would apply also to the Full Information case.\\
\indent The Principal-Agent problem under Hidden Action (or moral hazard) is inspired by the seminal paper of Holmstr\"{o}m and Milgrom \cite{MR882097} and is well treated for instance in \cite{MR2963805} and \cite{Williams}. In what follows we first consider a mean-variance version of such a model, using the stochastic maximum principle as our main tool. For reasons of tractability (and in line with the present literature) we will have to consult weak solutions of SDEs. We clarify this by referring to the model as \textit{Hidden Action in the weak formulation}.\\
\indent We will also consider a simpler model under Hidden Action in which the information set of the Principal is relaxed to a larger set. Such a relaxation does not necessarily imply Full Information and we refer to this model as \textit{Hidden Contract in the strong formulation}.

\subsection{Mean-Variance Hidden Action in the weak formulation}
\label{MVHiddenActionWF}
\indent Consider a Principal-Agent model where output $x(t)$ is modelled as a risky asset solving the SDE:
\begin{equation}
\left\{
    \begin{array}{l}
        dx(t) = \sigma(t,x(t)) dW_t,\\
        x(0) = 0,
    \end{array}
\right.
\end{equation}
Here $T > 0$ and $W_t$ is a 1-dimensional standard Brownian motion defined on the filtered probability space $(\Omega, \mathcal{F}, \mathbb{F}, \mathbb{P})$. For the diffusion we assume $\sigma > 0$ and $\mathbb{E}\left[ \int_0^t \sigma(t,x(t))^2 dt \right] < \infty$. The agent's level of effort is represented by a process $e(\cdot)$, taking values in some predefined subset $E \subseteq \R$ (typically $E = [0,\hat{e}]$ for some non-negative $\hat{e}$, or $E = \R$) and is required to belong to the set $\mathcal{E}[0,T]$, where
\begin{equation*}
    \mathcal{E}[0,T] := \{ e: [0,T] \times \Omega \rightarrow E; \hspace{0.2cm} e \text{ is } \mathbb{F} \text{-adapted} \}
\end{equation*}
We consider the case of Hidden Actions meaning that the principal cannot observe $e(\cdot)$. Output, however, is public information and observed by both the principal and the agent. Before the period starts the principal specifies an $\mathbb{F}_x$-adapted cash-flow $s(\cdot)$ (typically non-negative) for all $t \in [0,T]$, which compensates the agent for costly effort in managing $x(\cdot)$. Just as for the effort we assume $s(t) \in S$ for all $t \in [0,T]$ and some subset $S \subseteq \R$ and require $s(\cdot) \in \mathcal{S}[0,T]$, where
\begin{equation*}
    \mathcal{S}[0,T] := \{ s: [0,T] \times \Omega \rightarrow S;  \hspace{0.2cm} s \text{ is } \mathbb{F}_x \text{-adapted} \}.
\end{equation*}
The principal is not constrained by any means and can commit to any such process $s(\cdot) \in \mathcal{S}[0,T]$.\\
\indent In this model we consider cost functionals $\mathcal{J}_P$ and $\mathcal{J}_A$ of the principal and the agent respectively of the following form:
\begin{equation}
    \mathcal{J}_A (e(\cdot); s) := \E \left[ \int_0^T u (t,x(t),e(t),s(t)) dt + v(x(t)) \right]
    \label{TheProblemAgentCost}
\end{equation}
and
\begin{equation}
    \begin{array}{l}
        \vspace{0.2cm}\displaystyle \mathcal{J}_P(s(\cdot)) := -\E \left[ \int_0^T U (t, x(t), s(t)) dt +V(x(T)) \right]\\
        \hspace{1cm}\displaystyle+ \frac{r}{2} \text{Var} \left( \int_0^T \Phi (t, x(t), s(t)) dt + \Psi(x(T)) \right),
    \end{array}
\end{equation}
for some given risk aversion $r > 0$. The agent will accept $s(\cdot)$ and start working for the principal only if the participation constraint
\begin{equation}
    \mathcal{J}_A (\bar{e}(\cdot );s) \leq W_0,
    \label{TheProblemPartConst}
\end{equation}
is fulfilled by $s$ for some, typically negative, constant $W_0$. We assume \textit{incentive compatibility}, meaning that the agent will act as to optimize $\mathcal{J}_A$ in response to any given $s(\cdot)$. The principal's problem is to minimize $\mathcal{J}_P$ under the participation constraint and incentive compatibility.\\
\indent A direct approach to the Principal-Agent problem as described above is, however, not mathematically tractable. Therefore, in line with \cite{MR2963805} and \cite{Williams} we make the problem tractable by using the concept of weak solutions of SDEs. That is, rather than having a model in which the agent controls $x(\cdot)$ itself we consider the \textit{density of output}, $\Gamma^e(\cdot)$, as the controlled dynamics where
\begin{equation}
\left\{
    \begin{array}{l}
        d\Gamma^e(t) = \Gamma^e(t) f(t, x(t), e(t)) \sigma^{-1}(t, x(t)) dW_t,\\
        \Gamma^e(0) = 1,
    \end{array}
    \label{HiddenActionAgentDyn}
\right.
\end{equation}
for a given function $f$ describing production rate and satisfying Assumption 1 in Section \ref{Preliminaries}. Note that $\Gamma^e = \mathcal{E}(Y)$ where $Y(t) = \int_0^t f(\tau,x(\tau),e(\tau))\cdot \sigma^{-1}(\tau,x(\tau))dW(\tau)$ and $\mathcal{E}(\cdot)$ denotes the stochastic exponential. The key idea behind the weak formulation of the Hidden Action model, letting the agent control $\Gamma^e(t)$ rather than $x(t)$, is that it allows us to consider $(x)_T$ as a fixed but random realization (actually $\mathbb{F}_x = \mathbb{F}$ as a consequence of the regularity of $\sigma$). If $\Gamma^e(\cdot)$ is a martingale, which follows by assuming for instance the Novikov condition or the Bene\v{s} condition (see. \cite{MR1121940} p. 200), we have by Girsanov's theorem that the probability measure $d\mathbb{P}^e$ defined by
\begin{equation}
    \frac{d\mathbb{P}^e}{d\mathbb{P}} = \Gamma^e(T)
\end{equation}
makes the process $W^e(t)$ defined by
\begin{equation}
    dW^e_t = dW_t - f(t,x(t),e(t))\sigma^{-1}(t,x(t))dt
\end{equation}
a $\mathbb{P}^e$-Brownian motion. In particular
\begin{equation}
    dx(t) = f(t,x(t),e(t))dt + \sigma(t,x(t))dW^e_t
    \label{HiddenActionWeakDyn}
\end{equation}
and
\begin{equation}
    \begin{array}{l}
        \vspace{0.2cm}\displaystyle\mathcal{J}_A(e(\cdot);s) = \mathbb{E}^e \left[ \int_0^T u(t,x(t),e(t),s(t))dt + v(x(T)) \right] =\\
        \hspace{0.3cm}\displaystyle= \mathbb{E} \left[ \int_0^T \Gamma^e(t)u(t,x(t),e(t),s(t))dt + \Gamma^e(T)v(x(T)) \right].
    \end{array}
    \label{HiddenActionAgentCost}
\end{equation}
We think of the Principal-Agent problem as divided into two coupled problems; \textit{The Agent's problem} and \textit{The Principal's problem}.\vspace{0.2cm}\\
\noindent \textbf{The Agent's problem (weak formulation):}
Given any $s(\cdot) \in \mathcal{S}[0,T]$ (that we assume fulfills the participation constraint) the Agent's problem is to find a process $\bar{e}(\cdot) \in \mathcal{E}[0,T]$ such that the cost functional
\begin{equation*}
    \mathcal{J}_A (\bar{e}(\cdot );s) = \mathbb{E} \left[ \int_0^T \Gamma^e(t)u(t,x(t),e(t),s(t))dt + \Gamma^e(T)v(x(T)) \right],
\end{equation*}
is minimized, subject to the dynamics in (\ref{HiddenActionAgentDyn}).\vspace{0.2cm}\\
\noindent \textbf{The Principal's problem (strong formulation):}
Given that the Agent's problem has an optimal solution $\bar{e}(\cdot)$ in the weak formulation the Principal's problem is to find a process $\bar{s}(\cdot) \in \mathcal{S}[0,T]$,
such that the cost functional
\begin{equation*}
    \begin{array}{l}
        \vspace{0.2cm}\displaystyle \mathcal{J}_P(\bar{s}(\cdot)) := -\E \left[ \int_0^T U (t, x(t), s(t)) dt +V(x(T)) \right]\\
        \hspace{1cm}\displaystyle+ \frac{r}{2} \text{Var} \left( \int_0^T \Phi (t, x(t), s(t)) dt + \Psi(x(T)) \right),
    \end{array}
\end{equation*}
is minimized and
\begin{equation*}
    \mathcal{J}_A (\bar{e}(\cdot );\bar{s}) = \mathbb{E} \left[ \int_0^T u(t,x(t),e(t),s(t))dt + v(x(T)) \right] \leq W_0,
\end{equation*}
subject to the dynamics
\begin{equation*}
\left\{
    \begin{array}{l}
    dx(t) = \sigma (t, x(t)) dW_t, \hspace{0.5cm} t\in (0,T], \\
    x(0) = 0.
    \end{array}
\right.
\end{equation*}
\begin{remark}
    Here we have chosen to formulate the Principal's problem in the strong form rather than in the weak form, which seems to be most common in the literature. However, as pointed out in \cite{MR2963805}, because of adaptiveness this approach can be problematic in certain models. This is a fact that one should be aware of.
\end{remark}
\noindent In this context the following definition is natural.
\begin{definition}
An \textit{optimal contract} is a pair $( \bar{e}(\cdot), \bar{s}(\cdot) ) \in \mathcal{E}[0,T] \times \mathcal{S}[0,T]$ obtained by sequentially solving first the Agent's- and then the Principal's problem.
\end{definition}
\noindent In game theoretic terminology an optimal contract can thus be thought of as a Stackelberg-equilibrium in a two-player non-zero-sum game.\\
\indent It is important to note that even though the principal cannot observe the agent's effort, he/she can still offer the agent a contract by suggesting a choice of effort $e(\cdot)$ and a compensation $s(\cdot)$. By incentive compatibility, however, the principal knows that the agent only will follow such a contract if the suggested effort solves the agent's problem. To find the optimal effort, $\bar{e}(\cdot)$, the principal must have information of the agent's preferences, i.e. the functions $u$ and $v$. The realism of such an assumption is indeed questionable but nevertheless necessary in our formulation due to the participation constraint. In order to make the intuition clear and to avoid any confusion we adopt the convention that the principal has full information of the agent's preferences $u$ and $v$. This gives a tractable way of thinking of how actual contracting is realized.\\
\indent Thus, the principal is able to predict the optimal effort $\bar{e}(\cdot)$ of the agent's problem and thereby suggest an optimal contract $(\bar{e}(\cdot), \bar{s}(\cdot))$, if it exists.\\
\indent The idea is to apply the methods from Section \ref{Preliminaries} to characterize optimal contracts in the general Principal-Agent model presented above. However, since the control variable $e$ figures in the diffusion of (\ref{HiddenActionAgentDyn}) we require the following convexity assumption in order to avoid a second order adjoint process in the maximum principle:\vspace{0.2cm}\\
\noindent \textbf{Assumtion 7:} The set $E \subset \R$ is convex.\vspace{0.2cm}\\
The Agent's Hamiltonian in the weak formulation is
\begin{equation}
    \mathcal{H}_A(t,x,\Gamma^e,e,p,q,s) := q \cdot \Gamma^e \cdot \frac{f(t, x, e)}{\sigma(t, x)} - \Gamma^e \cdot u(t,x,e,s),
    \label{HiddenActionAgHam}
\end{equation}
and by Theorem \ref{PrelThmSMP1} any optimal control $\bar{e}(t)$ solving the Agent's problem must maximise $\mathcal{H}_A$ pointwisely. The pair $(p(\cdot), q(\cdot))$ solves the Agent's adjoint BSDE:
\begin{equation}
\left\{
    \begin{array}{l}
        \vspace{0.1cm}dp(t) =  -\left\{ q(t) \cdot \frac{f(t, x(t), \bar{e}(t))}{\sigma(t, x(t))} - u(t,x(t),\bar{e}(t),s(t)) \right\}dt + q(t) dW_t,\\
        p(T) = -v_x(x(T))
    \end{array}
    \right.
    \label{HiddenActionAdjBSDE}
\end{equation}
If $f$ and $u$ both are differentiable in the $e$ variable and we assume that $\bar{e}(\cdot) \in \text{int}(E)$, maximizing $\mathcal{H}_A$ translates into the first order condition
\begin{equation}
    q(t) = \sigma(t,x(t)) \cdot \frac{u_e(t,x(t),\bar{e}(t), s(t))}{f_e(t,x(t),\bar{e}(t))},
    \label{HiddenActionCharofe}
\end{equation}
which is in agreement with \cite{Williams}. Before proceeding to the Principal's problem we assume solvability of $\bar{e}$ in (\ref{HiddenActionCharofe}) and we write
\begin{equation*}
    \bar{e}(t) = e^*(t, \bar{x}(t), q(t), s(t)),
\end{equation*}
where $e^* : \R_+ \times \R^4 \rightarrow \R$ is a function having sufficient regularity to allow for the existence of a unique solution to the FBSDE (\ref{TheProblemFBSDE}) below. Based on the information given by $e^*$ the principal wishes to minimize the cost $\mathcal{J}_P$ by selecting a process $s(\cdot)$ respecting (\ref{TheProblemPartConst}). The dynamics of the corresponding control problem is, in contrast to the SDE of the agent's problem, a FBSDE built up by the output SDE coupled to the agent's adjoint BSDE. More precisely:
\begin{equation}
    \left\{
    \begin{array}{l}
        dx(t) =  \sigma (t, x(t)) dW_t,\vspace{0.1cm}\\
        dp(t) = -\left\{ q(t) \cdot \frac{f(t, x(t), e^*(t, x(t), q(t), s(t)))}{\sigma(t, x(t))} - u(t,x(t),e^*(t, x(t), q(t), s(t)),s(t)) \right\}dt\vspace{0.1cm}\\
        \hspace{1.3cm}+ q(t) dW_t,\vspace{0.1cm}\\
        \bar{x}(0) = 0, \hspace{0.2cm} p(T) = - v_x(x(T)).
    \end{array}
    \right.
    \label{TheProblemFBSDE}
\end{equation}
\noindent In order to characterize cash-flow optimality in the Principal's problem we apply Theorem \ref{MVSMPMV1}. The Hamiltonian reads
\begin{equation}
    \begin{array}{l}
    \vspace{0.1cm}\mathcal{H}_P(t,x,q,s,R,P_1, P_2,Q_1, Q_2, \lambda_P,\lambda_A) :=\\
    \vspace{0.1cm}R \cdot \left\{ -q \cdot \frac{f(t, x, e^*(t, x, q, s))}{\sigma(t, x)} + u(t,x,e^*(t, x, q, s),s) \right\} + P_2 \cdot (\Phi(t,x,s) +\frac{\sigma^2(t,x)}{2} \Psi''(x))\\
    + Q_1 \cdot \sigma (t, x) + Q_2 \cdot \sigma (t, x) \Psi'(x) - \lambda_A \cdot u(t,x,e^*(t, x, q, s),s) + \lambda_P \cdot \mathcal{U} (t,x,s),
    \end{array}
    \label{HiddenActionPrincHam}
\end{equation}
and for any optimal 4-tuple $(\bar{x}(\cdot), \bar{p}(\cdot), \bar{q}(\cdot), \bar{s}(\cdot))$ we have the existence of Lagrange multipliers $\lambda_A, \lambda_P \in \R$ satisfying the conditions in Theorem (\ref{MVSMPMV1}). The adjoint processes $(R(\cdot), P_1(\cdot), Q_1(\cdot), P_2(\cdot), Q_2(\cdot))$ solve the FBSDE (\ref{MVCor1FBSDE}), in which case
\begin{equation*}
        \bar{s}(t)= \smash{\displaystyle\arg\max_{s \in S}}\hspace{0.1cm} \mathcal{H}_P(t, \bar{x}(t),\bar{p}(t),\bar{q}(t),s,R(t),P_1(t), Q_1(t), P_2(t), Q_2(t),\lambda_P,\lambda_A).
\end{equation*}
Before stating the full characterization of optimal contracts in the Mean-Variance Principal-Agent problem under Hidden Action we introduce the following technical assumption:\\
\newline
\textbf{(PA1).} All functions involved in the Agent's problem satisfy Assumption 1 from Section \ref{Preliminaries} and the density of output is a martingale. The functions defining the Principal's problem (including composition with the map $e^*$) satisfy Assumptions 2-6, also from Section \ref{Preliminaries}, and $\Psi$ is three times differentiable.\\
\begin{theorem}
    Let the statements in (PA1) and Assumption 7 hold and consider the Mean-Variance Principal-Agent problem under Hidden Actions with risk aversion $r > 0$ and participation constraint defined by $W_0 < 0$. Then, if $(\bar{e}(\cdot), \bar{s}(\cdot))$ is an optimal contract there exist numbers $\lambda_A, \lambda_P \in \R$ such that
        \begin{equation*}
            \lambda_P \geq 0, \hspace{0.5cm} \lambda_A^2 + \lambda_P^2 = 1,
        \end{equation*}
    a pair $(p(\cdot), q(\cdot))\in L_{\mathcal{F}}^2(0,T;\R) \times (L_{\mathcal{F}}^2(0,T;\R))$ solving the SDE in (\ref{HiddenActionAdjBSDE}) and a quintuple $(R(\cdot), P_1(\cdot), P_2(\cdot),Q_1(\cdot), Q_2(\cdot)) \in L_{\mathcal{F}}^2(\Omega; C([0,T]; \mathbb{R})) \times L_{\mathcal{F}}^2(\Omega; C([0,T]; \mathbb{R})) \times L_{\mathcal{F}}^2(0,T;\mathbb{R})$ solving the adjoint FBSDE (\ref{MVCor1FBSDE}) defined by (\ref{TheProblemFBSDE}) such that, sequentially,
        \begin{equation*}
            \bar{e}(t) = \smash{\displaystyle\arg\max_{e \in E}}\hspace{0.1cm}\mathcal{H}_A(t,\bar{x}(t),\Gamma^e(t),e,q(t),s(t)),
        \end{equation*}
    and
        \begin{equation*}
            \bar{s}(t)= \smash{\displaystyle\arg\max_{s \in S}}\hspace{0.1cm} \mathcal{H}_P(t,\bar{x}(t),\bar{q}(t),s,R(t),P_1(t), Q_1(t),P_2(t), Q_2(t),\lambda_P,\lambda_A),
        \end{equation*}
    with Hamiltonians $\mathcal{H}_A$ and $\mathcal{H}_P$ as in (\ref{HiddenActionAgHam}) and (\ref{HiddenActionPrincHam}) respectively.
    \label{HiddenActionMainThm}
\end{theorem}

\subsection{Hidden Contract in the strong formulation}
\label{MVHiddenContract}
We are now going to study a different type of Mean-Variance Principal-Agent problems called \textit{Hidden Contract} models (introduced in \cite{PA}). Comparing to the Hidden Action model in Section \ref{MVHiddenActionWF} the Hidden Contracts differ in two key aspects. First we relax the information set of the Principal from $\mathbb{F}_x$ to the full filtration generated by the Brownian motion. Secondly we treat the process $s(\cdot)$ as hidden, meaning that the Agent reacts to the provided cash-flow given as an $\mathbb{F}$-adapted process, without being aware of the underlying dependence of output. This explains the name Hidden Contract.\\
\indent
\indent The fact that the underlying mathematical structure of $s(\cdot)$ is unknown to the Agent in the Hidden Contract model motivates the relevance of a Mean-Variance framework by an extended participation constraint (compared to (\ref{TheProblemPartConst})). By requiring an upper bound for the variance of for instance the expected accumulated wealth provided by $s(\cdot)$ the Agent can protect him/her-self from undesirable high levels of risk. The setup goes as follows:\\
\indent Consider a Principal-Agent model in which output $x(t)$ is modelled as a risky asset solving the SDE
\begin{equation}
\left\{
    \begin{array}{l}
    dx(t) = f(t, x(t), e(t)) dt + \sigma (t, x(t)) dW_t, \hspace{0.5cm} t\in (0,T], \\
    x(0) = 0.
    \end{array}
\right.
\end{equation}
Here $T > 0$ and $W_t$ is a 1-dimensional standard Brownian motion defined on the filtered probability space $(\Omega, \mathcal{F}, \mathbb{F}, \mathbb{P})$. The functions $f$ and $\sigma$ represent production rate and volatility respectively, and we assume both of them to satisfy Assumption 1 from the Section \ref{Preliminaries}. Just as for the Hidden Action case we require any admissible effort process $e(\cdot)$ to be in $\mathcal{E}[0,T]$. For the admissible cash-flows, however, we enlarge $\mathcal{S}[0,T]$ (due to the extended flow of information to the Principal) to
\begin{equation*}
    \mathcal{S}[0,T] := \{ s: [0,T] \times \Omega \rightarrow S;  \hspace{0.2cm} s \text{ is } \mathbb{F} \text{-adapted} \}.
\end{equation*}
We consider the cost functionals
\begin{equation}
    \mathcal{J}_A (e(\cdot); s) := \E \left[ \int_0^T u (t,x(t),e(t),s(t)) dt + v(x(T)) \right],
    \label{MVHiddenContractAgCost}
\end{equation}
and
\begin{equation}
    \mathcal{J}_P(s(\cdot)) := \E \left[ \int_0^T \mathcal{U} (t,x(t),s(t)) dt + \mathcal{V}(x(T)) \right],
    \label{MVHiddenContractPrincCost}
\end{equation}
and the participation constraint:
\begin{equation}
    \left\{
    \begin{array}{l}
        \displaystyle\vspace{0.1cm}\mathcal{J}_{A}(\bar{e}(\cdot);s) := \E \left[ \int_0^T u (t,x(t),\bar{e}(t),s(t)) dt + v(x(T)) \right] \leq W_0,\\
        \displaystyle\mathcal{I}_{A}(\bar{e}(\cdot);s) := \text{Var}\left( \int_0^T \phi (t,x(t),\bar{e}(t),s(t)) dt + \psi(x(T)) \right) \leq R_0.
    \end{array}
    \right.
    \label{MVHiddenContractPartConst}
\end{equation}
Just as for the Hidden Action case in Section \ref{MVHiddenActionWF} we consider the Agent's- and the Principal's problem sequentially. The precise statements are:\\
\newline
\noindent \textbf{The Agent's Problem.} Given any $s(\cdot) \in \mathcal{S}[0,T]$ (fulfilling the participation constraint) the Agent's problem is to find a process $\bar{e}(\cdot) \in \mathcal{E}[0,T]$ minimizing (\ref{MVHiddenContractAgCost}).\\
\newline
\noindent \textbf{The Principal's Problem.} Given that the Agent's problem has an optimal solution $\bar{e}(\cdot)$ the Principal's problem is to find a process $\bar{s}(\cdot) \in \mathcal{S}[0,T]$ minimizing the cost functional (\ref{MVHiddenContractPrincCost}) subject to the participation constraint (\ref{MVHiddenContractPartConst}).\\
\newline
The mathematical virtue of Hidden Contracts is the possibility of working solely in the strong formulation. For the Agent's problem we are facing the Hamiltonian
\begin{equation}
    \mathcal{H}_A (t, x, e, p, q, s) := p \cdot f(t, x, e) + q \cdot \sigma(t, x) - u(t, x, e, s).
    \label{MVHiddenContractAgHam}
\end{equation}
Therefore, by Theorem \ref{PrelThmSMP1} we have for any optimal pair $(\bar{x}(\cdot), \bar{e}(\cdot))$ the existence of adjoint processes $(p(\cdot), q(\cdot))$ solving the BSDE:
\begin{equation}
\left\{
    \begin{array}{l}
        \vspace{0.1cm}dp(t) = - \left\{ f_x(t, \bar{x}(t), \bar{e}(t)) p(t) + \sigma_x(t, \bar{x}(t)) q(t) - u_x (t, \bar{x}(t), \bar{e}(t)) \right\} dt + q(t) dW_t, \\
        p(T) = - v_x(\bar{x}(T)),
    \end{array}
    \label{MVHiddenContractAgentBSDE}
\right.
\end{equation}
and the characterization
\begin{equation}
    \bar{e}(t)= \smash{\displaystyle\arg\max_{e \in E}}\hspace{0.1cm}\mathcal{H}_A(t, \bar{x}(t), e, p(t), q(t), s(t)),
    \label{MVHiddenContractAgentHamMax}
\end{equation}
for a.e. $t \in [0,T]$ and $\Prob$-a.s.\\
\indent As in the Hidden Contract case we proceed into the Principal's problem by assuming the existence of a function $e^*$ such that $\bar{e}_t = e^*(t, x(t), p(t), q(t), s(t))$ (having sufficient regularity to allow for existence and uniqueness of a solution to (\ref{MVHiddenContractFBSDE})). The Principal is facing the problem of minimizing $\mathcal{J}_P$ subject to (\ref{MVHiddenContractPartConst}) by controlling the following FBSDE:
\begin{equation}
    \left\{
    \begin{array}{l}
        \vspace{0.1cm}d\bar{x}_t = f(t, \bar{x}(t), e^*(t, \bar{x}(t), p(t), q(t), s(t))) dt + \sigma (t, \bar{x}(t)) dW_t,\\
        \vspace{0.1cm}dp(t) = - \left\{ f_x(t, \bar{x}(t), e^*(t, \bar{x}(t), p(t), q(t), s(t))) p(t) + \sigma_x(t, \bar{x}(t)) q(t) \right.\\
        \vspace{0.1cm}\hspace{4cm}\left. - u_x (t, \bar{x}(t), e^*(t, \bar{x}(t), p(t), q(t), s(t))) \right\} dt + q(t) dW_t,\\
        \bar{x}(0) = 0, \hspace{0.2cm} p(T) = - v_x(\bar{x}(T)).
    \end{array}
    \right.
    \label{MVHiddenContractFBSDE}
\end{equation}
We now apply Theorem \ref{MVSMPMV2} in order to characterize optimal cash-flows in the Principal's problem. The associated Hamiltonian is
\begin{equation}
    \begin{array}{l}
    \vspace{0.1cm}\mathcal{H}_P(t,x,p,q,s,R,P_1,P_2,Q_1,Q_2,\lambda_E, \lambda_V,\lambda_P) =\\
    \vspace{0.1cm}R \cdot (f_x(t,x,e(t,x,p,q,s)) p + \sigma_x(t,x) q - u_x(t, x, e(t,x,p,q,s), s))\\
    \vspace{0.1cm}+ P_1 \cdot f(t,x,e(t,x,p,q,s)) + P_2 \cdot \bigg\{ \phi(t, x, e(t,x,p,q,s), s)\\
    \vspace{0.1cm}+ f(t,x,e(t,x,p,q,s))\psi'(x) + \frac{\sigma^2(t,x)}{2} \psi''(x) \bigg\} + Q_1 \cdot \sigma(t,x) + Q_2 \cdot \sigma(t,x) \psi'(x)\\
    - \lambda_E \cdot u(t, x, e(t,x,p,q,s), s) - \lambda_P \cdot \mathcal{U}(t,x,s).
    \end{array}
    \label{HiddenContractPrincHam}
\end{equation}
For any optimal 4-tuple $(\bar{x}(\cdot), \bar{p}(\cdot),\bar{q}(\cdot), \bar{s}(\cdot))$ of the Principal's problem we have the existence of Lagrange multipliers $\lambda_E, \lambda_V, \lambda_P \in \R$ satisfying either of the conditions (i)-(v) in Section \ref{Preliminaries}, with $\lambda_P \geq 0$ and
\begin{equation*}
    \lambda_E^2 + \lambda_V^2 + \lambda_P^2 = 1,
\end{equation*}
and a triple of adjoint processes $(R(\cdot), P(\cdot), Q(\cdot))$ solving the FBSDE (\ref{MVCor2FBSDE}) so that
\begin{equation*}
    \vspace{0.1cm}\bar{s}(t) = \smash{\displaystyle\arg\max_{s \in S}}\hspace{0.1cm} \mathcal{H}_P(t,x(t),p(t),q(t),s,R(t),P_1(t),P_2(t),Q_1(t),Q_2(t),\lambda_E,\lambda_P).
\end{equation*}
For the full characterization of optimality we require the following technical assumption:\\
\newline
\textbf{(PA2).} All functions involved in the Agent's problem satisfy the Assumption 1 from Section \ref{Preliminaries}. The functions defining the Principal's problem (including composition with the map $e^*$) satisfy the Assumptions 2-6, also from Section \ref{Preliminaries}, and $\psi$ is three times differentiable.\\
\begin{theorem}
    Let the statements in (PA2) hold and consider the Mean-Variance Principal-Agent problem under Hidden Contract with participation constraints defined by the given parameters $W_0 < 0$ and $R_0 > 0$. Then, if $(\bar{e}(\cdot), \bar{s}(\cdot))$ is an optimal contract there exist numbers $\lambda_E, \lambda_V, \lambda_P \in \R$ such that
        \begin{equation*}
            \lambda_P \geq 0, \hspace{0.5cm} \lambda_E^2 + \lambda_V^2 + \lambda_P^2 = 1,
        \end{equation*}
    a pair $(p(\cdot), q(\cdot))\in L_{\mathcal{F}}^2(0,T;\R) \times (L_{\mathcal{F}}^2(0,T;\R))$ solving the BSDE in (\ref{MVHiddenContractAgentBSDE}) and a quintuple $(R(\cdot), P_1(\cdot), P_2(\cdot),Q_1(\cdot), Q_2(\cdot)) \in L_{\mathcal{F}}^2(\Omega; C([0,T]; \mathbb{R})) \times L_{\mathcal{F}}^2(\Omega; C([0,T]; \mathbb{R})) \times L_{\mathcal{F}}^2(0,T;\mathbb{R})$ solving the adjoint FBSDE (\ref{MVCor2FBSDE}) defined by (\ref{MVHiddenContractFBSDE}) such that, sequentially,
        \begin{equation*}
            \bar{e}(t) = \smash{\displaystyle\arg\max_{e \in E}}\hspace{0.1cm}\mathcal{H}_A(t, \bar{x}(t), e, p(t), q(t), s(t)),
        \end{equation*}
    and
        \begin{equation*}
            \vspace{0.1cm}\bar{s}(t) = \smash{\displaystyle\arg\max_{s \in S}}\hspace{0.1cm} \mathcal{H}_P(t,\bar{x}(t),\bar{p}(t),\bar{q}(t),s,R(t),P_1(t),P_2(t),Q_1(t),Q_2(t),\lambda_E,\lambda_V, \lambda_P).
        \end{equation*}
    with Hamiltonians $\mathcal{H}_A$ and $\mathcal{H}_P$ as in (\ref{MVHiddenContractAgHam}) and (\ref{HiddenContractPrincHam}) respectively.
    \label{HiddenContractMainThm}
\end{theorem}

\section{A Solved Example in the Case of Hidden Contracts}
\label{HiddenContractLQ}
\begin{figure}[ht!]
\includegraphics[width=1.0\textwidth]{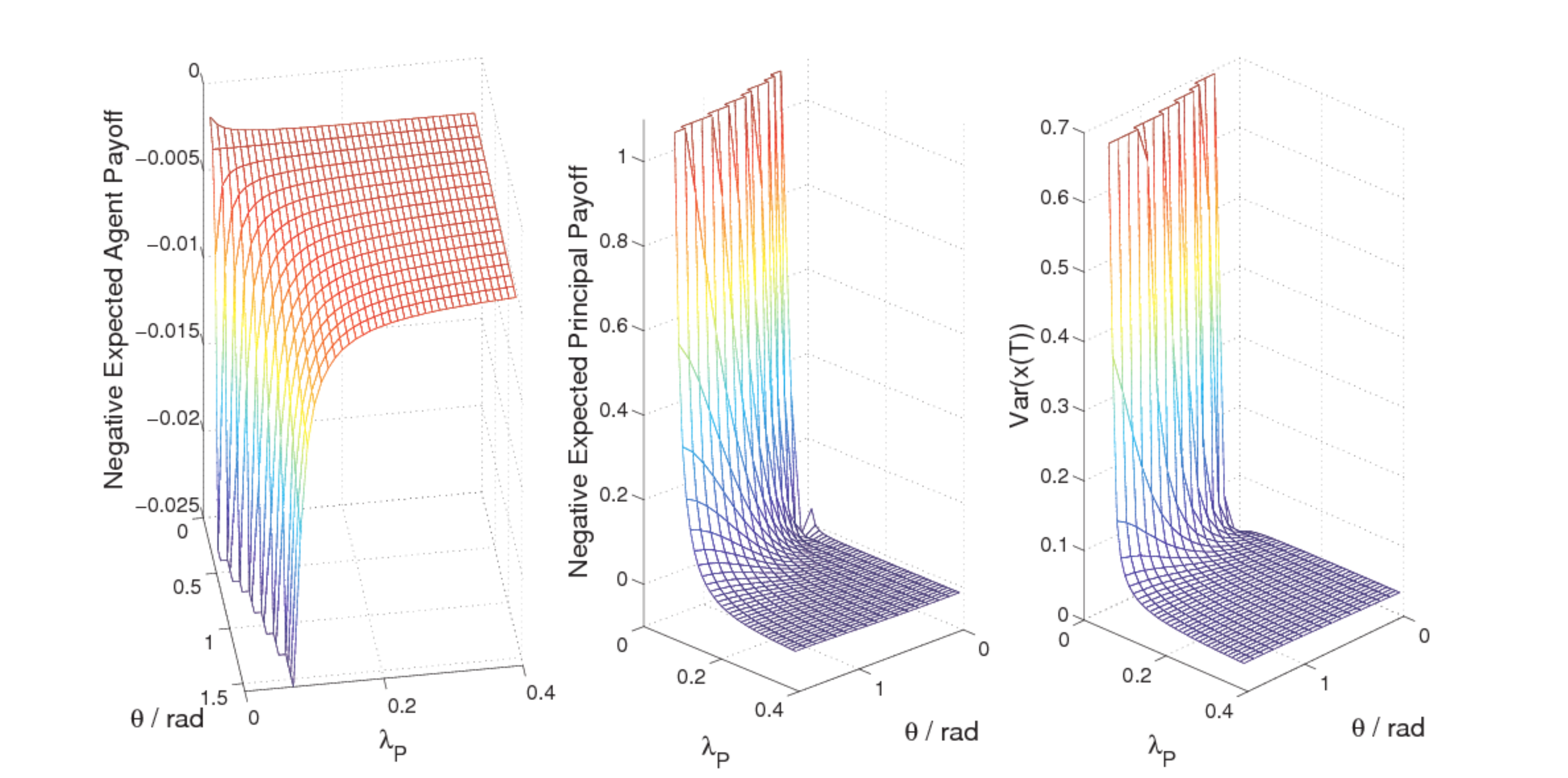}
  \caption{Monte-Carlo simulations of $\mathcal{J}_A(\bar{e}(\cdot);\bar{s})$, $\mathcal{J}_P(\bar{s}(\cdot))$ and $\text{Var}(x(T))$ as functions of $\lambda_P$ and $\theta$ (relating to $\lambda_E$ and $\lambda_V$ via case (iv)) based on $10^6$ sample paths at each point. Parameter values: $a = b = \sigma = 1, \alpha = 0.2, \beta = 1, T = 0.03$.}
  \label{HiddenContractLQExpPayoffFig}
\end{figure}

We now illustrate the method of Section \ref{TheProblem} by considering a concrete example of Hidden Contract type. In order to find explicit solutions we choose a linear-quadratic setup. As a result we get optimal contracts adapted to the filtration generated by output.\\
\indent Consider the following dynamics of production,
\begin{equation*}
    \left\{
    \begin{array}{l}
        dx(t) = (ax(t) + be(t))dt + \sigma dW_t, \hspace{0.3cm} t\in (0,T],\\
        x(0) = 0, \hspace{0.3cm} a,b \in \R \text{ and } \sigma > 0,
    \end{array}
    \right.
\end{equation*}
and let the preferences of the agent and the principal be described by quadratic utility functions:
\begin{eqnarray}
    && \mathcal{J}_A (e(\cdot);s) := \E\left[ \int_0^T \frac{(s_t - e_t)^2}{2} dt - \alpha \cdot \frac{x(T)^2}{2} \right],\\
    && \mathcal{J}_P (s(\cdot)) := \E\left[\int_0^T \frac{s_t^2}{2} dt - \beta \cdot \frac{x(T)^2}{2}  \right].
    \label{HiddenContractLQPayoff}
\end{eqnarray}
Note that we are following the convention of Section \ref{TheProblem} to consider cost- rather than payoff-functionals. Thus, the Agent's utility function should be interpreted as a desire to maintain a level of effort close to the compensation given by the cash-flow.
We think of the parameters $\alpha>0$ and $\beta>0$ as bonus factors of total production at time $T$. For the participation constraint we require any admissible cash-flow $s(t)$ to satisfy the following:
\begin{equation}
\left\{
    \begin{array}{l}
        \mathcal{J}_A(\bar{e}(\cdot); s) \leq W_0,\\
        \text{Var} \left( x(T) \right) < R_0,
    \end{array}
\right.
    \label{HiddenContractLQConstr}
\end{equation}
where $W_0 < 0$, $R_0 > 0$ and $\bar{e}(\cdot)$ denotes the optimal effort policy of the agent given $s(\cdot)$.\\
\indent Assume that the principal offers the agent $s(\cdot)$ over the period $0 \leq t \leq T$. The Hamiltonian function of the agent is
\begin{equation*}
    \mathcal{H}_A(x,e,p,q,s) := p \cdot (ax + be) + q \cdot \sigma - \frac{(s-e)^2}{2},
\end{equation*}
so
\begin{equation}
    \frac{\partial \mathcal{H}_A}{\partial e} = bp + s - e = 0 \hspace{0.5cm} \text{and} \hspace{0.5cm} \bar{e}(t)= bp(t) + s(t),
    \label{HiddenContractLQOptEff}
\end{equation}
where the pair $(p,q)$ solves the adjoint equation
\begin{equation*}
    \left\{
    \begin{array}{l}
        \vspace{0.1cm}dp(t) = -ap(t) dt + q(t) dW_t,\\
        p(T) = \alpha x(T).
    \end{array}
    \right.
\end{equation*}
\noindent Turning to the principal's problem we want to control the FBSDE
\begin{equation}
    \left\{
    \begin{array}{l}
        dx(t) = (ax(t) + b^2p(t) + bs(t))dt + \sigma dW_t,\\
        dp(t) = -ap(t)dt + q(t)dW_t,\\
        x(0)=0, p(T)=\alpha x(T),
    \end{array}
    \label{HiddenContractLQFBSDE1}
    \right.
\end{equation}
optimally with respect to the cost function (\ref{HiddenContractLQPayoff}) and the participation constraint (\ref{HiddenContractLQConstr}). The Principal's Hamiltonian is
\begin{equation}
    \begin{array}{l}
        \vspace{0.1cm}\mathcal{H}_P (x,p,s,R,P_1,P_2,Q_1,Q_2,\lambda_E,\lambda_P) :=\\
        \hspace{0.2cm}\vspace{0.1cm}-ap \cdot R + (ax + b^2p + bs) \cdot P_1 + (s + ax + b^2p + bs) \cdot P_2 + \sigma \cdot (Q_1 + Q_2)\\
        \hspace{0.2cm}- \lambda_E \cdot \frac{b^2p^2}{2} - \lambda_P \cdot \frac{s^2}{2}
    \end{array}
\end{equation}
so
\begin{equation*}
    \frac{\partial \mathcal{H}_P}{\partial s} = bP_1 + (1+b)P_2 - \lambda_Ps \hspace{0.5cm} \text{and} \hspace{0.5cm} \bar{s}(t) = \frac{bP_1(t) + (1+b)P_2(t)}{\lambda_P},
\end{equation*}
where the quintuple $(R(t), P_1(t), P_2(t), Q_1(t), Q_2(t))$ solves the adjoint FBSDE:
\begin{equation}
\left\{
    \begin{array}{l}
        \vspace{0.1cm}dR(t) = (aR(t) - b^2(P_1(t) + P_2(t)) + \lambda_E b^2p(t))dt,\\
        \vspace{0.1cm}dP_1(t) = -a(P_1(t) + P_2(t))dt + Q_1(t)dW_t,\\
        \vspace{0.1cm}dP_2(t) = Q_2(t)dW_t,\\
        \vspace{0.1cm}R(0) = 0,\\
        P_1(T) = -\alpha R(T) + (\alpha\lambda_E +  \beta\lambda_P)x(T), \hspace{0.1cm}P_2(T) = 2\lambda_V(\mathbb{E}[\eta(T)] - \eta(T)).
    \end{array}
\right.
\label{HiddenContractLQBSDE}
\end{equation}
In this case, however, the auxiliary process $\eta(t)$ is the same as the output $x(t)$ in which case $P_2(T) = 2\lambda_V(\mathbb{E}[x(T)] - x(T))$. To solve the BSDE in (\ref{HiddenContractLQBSDE}) we can make a general linear ansatz:
\begin{equation}
    \left\{
    \begin{array}{l}
        \vspace{0.1cm}p(t) = A_{11}(t)x(t) + B_{11}(t)R(t) + A_{21}(t)\E [x(t)] + B_{21}(t)\E [R(t)],\\
        \vspace{0.1cm}P_1(t) = A_{12}(t)x(t) + B_{12}(t)R(t) + A_{22}(t)\E [x(t)] + B_{22}(t)\E [R(t)],\\
        P_2(t) = A_{13}(t)x(t) + B_{13}(t)R(t) + A_{23}(t)\E [x(t)] + B_{23}(t)\E [R(t)].\\
    \end{array}
    \right.
    \label{HiddenContractLQAnsatz}
\end{equation}
 Using the standard procedure with It\^{o}'s lemma it is elementary (but tedious) to derive a set of twelve coupled Riccati equations for the coefficients in (\ref{HiddenContractLQAnsatz}). A numerical example is presented in Figure \ref{HiddenContractLQriccatisystfig} below.
\begin{figure}[h!]
  \centering
\includegraphics[width=1.0\textwidth]{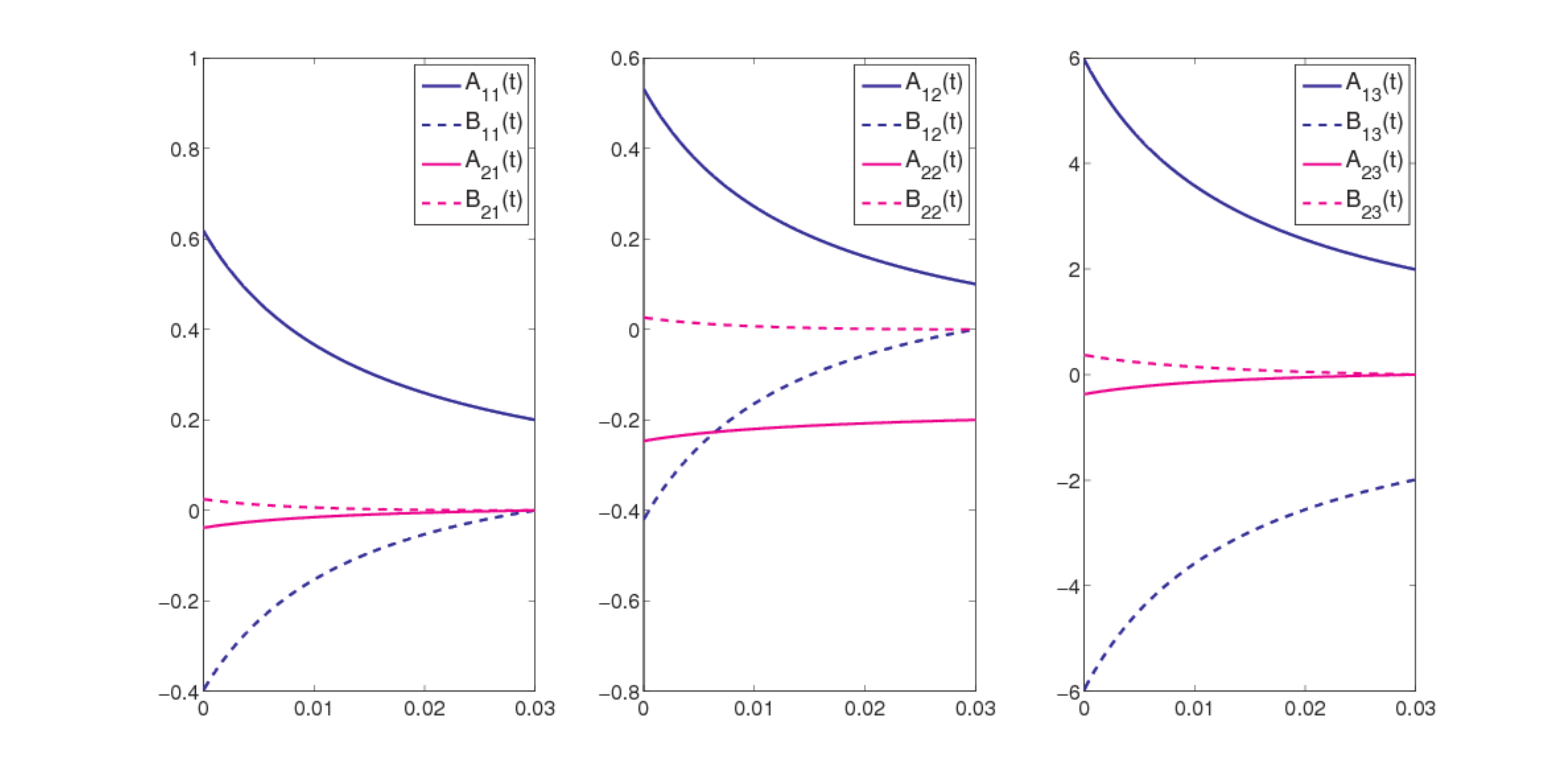}%
  \caption{Solution curves of (\ref{HiddenContractLQBSDE}) with parameter values chosen as: $a = b = \sigma = 1, \alpha = 0.2, \beta = 1, \lambda_P = 0.1, \theta = \pi/2, T = 0.03.$}
  \label{HiddenContractLQriccatisystfig}
\end{figure}
We get the unique semi-explicit solution to the optimal contract $\{ \bar{e}(t), \bar{s}(t) \}$, driven by the optimal dynamics $(\bar{x}(t), \bar{R}(t))$. What remains is to find a feasible tripple $(\lambda_E, \lambda_V, \lambda_P)$ so that the optimal contract fulfills the participation constraint in (\ref{HiddenContractLQConstr}). One way of finding such a triple is for instance by stochastic simulation of $(\bar{x}(t), \bar{R}(t))$ (e.g. a simple Euler-Maruyama scheme) and then estimate the payoff and the variance in (\ref{HiddenContractLQConstr}) by Monte-Carlo techniques for different values of $\lambda_P$. In Fig. \ref{HiddenContractLQExpPayoffFig} we have included the results of such a scheme corresponding to case (iv) of the transversality condition in Corollary \ref{MVSMPMV2}.
Note that $\bar{R}(t)$ satisfies the linear ODE
\begin{equation}
\left\{
\begin{array}{l}
    \frac{dR}{dt} + (b^2 B_{12} + b^2 B_{13}  - \lambda_E b^2 B_{11} - a) R(t) = (\lambda_E b^2 A_{11} - b^2 A_{12} - b^2 A_{13}) \bar{x}(t),\\
    R(0) = 0,
    \end{array}
    \right.
\end{equation}
so
\begin{equation*}
    \begin{array}{l}
    \vspace{0.1cm}R(t) =\\
     \displaystyle\frac{\int_0^t \text{exp} \left\{ \int_0^s b^2 B_{12} + b^2 B_{13}  - \lambda_E b^2 B_{11} - a \hspace{0.1cm} du \right\} \cdot (\lambda_E b^2 A_{11} - b^2 A_{12} - b^2 A_{13}) \bar{x} ds }{\text{exp} \left\{ \int_0^t b^2 B_{12} + b^2 B_{13}  - \lambda_E b^2 B_{11} - a \hspace{0.1cm} ds \right\}},
    \end{array}
\end{equation*}
and is by that $\mathbb{F}_x$-adapted. Therefore, in this model the optimal contract $\{ \bar{e}(t), \bar{s}(t) \}$ is $\mathbb{F}_x$-adapted and coincides with the corresponding strong solution to the Hidden Action problem, i.e. when the information set of the Principal is generated by output.\\

\end{document}